\documentclass{amsart}

\usepackage{graphicx}
\usepackage{amssymb}
\usepackage[text={6in,8.5in},centering]{geometry}
\usepackage{changes}
\newcommand{\R}{{\mathbb{R}}}
\newcommand{\Z}{{\mathbb{Z}}}
\newcommand{\N}{{\mathbb{N}}}

\newcommand{\ab}{{\mathbf{a}}}
\newcommand{\e}{{\mathbf{e}}}
\newcommand{\E}{{\mathbb{E}}}
\renewcommand{\P}{{\mathbb{P}}}
\newcommand{\X}{{\mathbf x}}
\newcommand{\y}{{\mathbf{y}}}
\newcommand{\z}{{\mathbf{z}}}
\renewcommand{\b}{{\mathbf{b}}}
\renewcommand{\S}{{\mathbb{S}}}
\newcommand{\Q}{{Q}}

\newtheorem{theorem}{Theorem}[section]

\newtheorem{remark}[theorem]{Remark}

\title[Ranges of random walk up to the time of exit]{Remarks on the range and multiple range of random walk up to the time of exit}

\author{Thomas Doehrman}
\address{Department of Mathematics\\
University of Arizona\\
621 N. Santa Rita Ave.\\
Tucson, AZ 85750, USA}
\email{{\tt thomasdoehrman@math.arizona.edu}}

\author{Sunder Sethuraman}
\address{Department of Mathematics\\
University of Arizona\\
621 N. Santa Rita Ave.\\
Tucson, AZ 85750, USA}
\email{{\tt sethuram@math.arizona.edu}}

\author{Shankar C. Venkataramani}
\address{Department of Mathematics\\
University of Arizona\\
621 N. Santa Rita Ave.\\
Tucson, AZ 85750, USA}
\email{{\tt shankar@math.arizona.edu}}
\date{}

\begin{document}

\begin{abstract}
We consider the scaling behavior of the range and $p$-multiple range, that is the number of points visited and the number of points visited exactly $p\geq 1$ times, of simple random walk on $\Z^d$, for dimensions $d\geq 2$, up to time of exit from a domain $D_N$ of the form $D_N = ND$ where $D\subset \R^d$, as $N\uparrow\infty$.  Recent papers have discussed connections of the range and related statistics with the Gaussian free field, identifying in particular that the distributional scaling limit for the range, in the case $D$ is a cube in $d\geq 3$, is proportional to the exit time of Brownian motion.   The purpose of this note is to give a concise, different argument that the scaled range and multiple range, in a general setting in $d\geq 2$, both weakly converge to proportional exit times of Brownian motion from $D$, and that the corresponding limit moments are `polyharmonic', solving a hierarchy of Poisson equations.    
\end{abstract}
\maketitle

{\scriptsize

}

\section{Introduction and results}
\label{results}
 
Let $\{X_n: n\geq 0\}$ be a simple random walk on $\Z^d$:  That is, 
$$\P\big(X_{n+1}=\X\pm \e_i|X_{n}=\X\big) = \frac{1}{2d},$$ 
for $1\leq i\leq d$ where $\big\{\e_i = (0,\ldots,0,1,0,\ldots, 0): 1\leq i\leq d\big\}$ is the standard basis of $\Z^d$.  The range of random walk up to time $n$, denoted by $\mathcal{R}_n$, is the number of distinct sites visited up to time $n$.
Correspondingly, for integers $p\geq 1$, the $p$-multiple range of random walk up to time $n$, denoted by $\mathcal{R}^{(p)}_n$, is the number of sites visited exactly $p$ times up to time $n$.  

To contrast with the range, the multiple range is a more delicate object.  For instance, whereas the set visited by the random walk up to time $n$ is connected and this range $\mathcal{R}_n$ is monotone in $n$, the analogous properties are not true for the $p$-multiple range set $\mathcal{R}^{(p)}_n$ that is visited exactly $p$ times.

\begin{figure}
  \includegraphics[trim=20 60 20 20, clip,width= \linewidth]{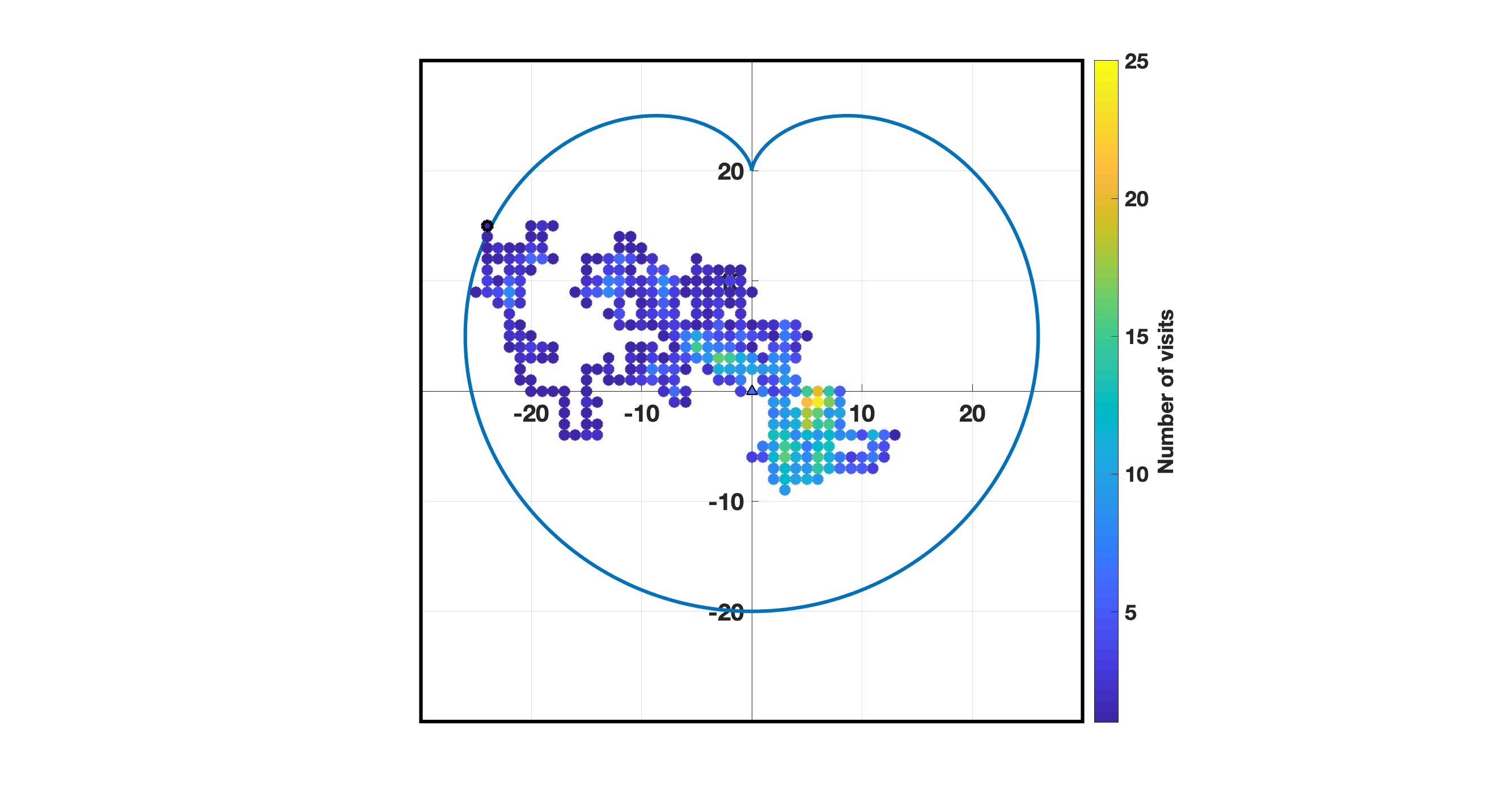}
  \caption{Range $\mathcal{R}_{\tau_N}$ and $p$-multirange $\mathcal{R}^{(p)}_{n}$ for a simple random walk on $\mathbb{Z}^2$, starting at the origin, up to the exit time from a cardioid $D$. For this realization, the exit time $n=936$ and $\mathcal{R}^{(p)}_n = 0$ unless $1\leq p\leq 24$.}
  \label{fig:multi}
\end{figure}

Although there is a large literature on the range and multiple range of random walk, often motivated by its applications and connections with other fields (cf. \cite{Asselah}, \cite{Bass}, \cite{Flatto}, \cite{Hamana}, \cite{Weiss}, \cite{Pitt},  and references therein),
the range and multiple range subject to constraints are less studied.  In particular, the range and multiple range up to the exit time from a domain have natural interpretations, which depend interestingly on the starting point of the random walk and the shape of the domain.

Recently, a body of work (see \cite{AB}, \cite{J1}, \cite{J2}) has considered extremes and local times of random walk in such settings and connections with the Gaussian free field.  We mention that \cite{J1}, \cite{J2} consider scaling limits of the number of `thick' points of continuous time random walk, starting at the origin, before its exit from a domain in $d\geq 2$. Thick points are those with visitation at least of order $a(\log N)^2$ in $d=2$ or of order $a \log N$ in $d \geq 3$, where $a\geq 0$ and $N$ is the length scale. In particular, when $a=0$ and the domain is a cube in $d\geq 3$, \cite[Theorem 1.5]{J2} shows that the scaled range $R_N/N^2$ converges weakly to a constant times the exit time of Brownian motion.

  In this context, the purpose of this note is to present a concise, different argument in a more general setting (Theorems \ref{thm1} and \ref{thm2}) that the distributional limits of the scaled range and multiple range up to the exit time from a domain in $d\geq 2$ connect to the exit times of Brownian motion and that their moments, as functions of the starting point, are `polyharmonic' in that they solve a hierarchy of Poisson PDEs.  In dimension $d=2$, the scalings are different, involving extra `$\log N$' factors, than in $d\geq 3$.  We also comment that the $p$-multiple range considers points with exactly $p$ visitations, where $p$ is fixed independent of $N$, a much different, and complementary object than the number of `thick' points considered previously.

Here, the proofs of Theorems \ref{thm1} and \ref{thm2} make use of the known scaling behavior of the unconstrained range and multiple range, the functional central limit theorem for random walks, and moment bounds that we provide in this note through simple discrete Fourier analysis.  Moreover, the argument shows that the randomness in the limits arises from the variability of the time to exit.  In particular, we observe that these moment bound derivations, in terms of a notion of `conductance', give a short way
to see the scalings in Theorems \ref{thm1} and \ref{thm2}.

We note that the phenomenon in $d\geq 2$ is different than in dimension $d=1$, considered before in \cite{locker} for the range, where the ordering of space forces a different type of limit (cf. Remark \ref{1D}).  The methods in \cite{locker} in $d=1$, make use of the reduced geometry, and do not carry over to higher dimensions $d\geq 2$.

\medskip

Let $N\geq 1$ be a scaling parameter, and let $D_N = ND$ where $D\subset \R^d$ is an open, bounded domain with some regularity.  To be definite, we will suppose that $D$ is a Lipschitz domain, although in $d=2$ it may be taken as a Jordan domain with a rectifiable boundary (which includes the case that $D$ is Lipschitz).
With respect to values $N\geq 1$, consider an array of simple random walks $\{X^{(N)}_n: n\geq 0\}$, starting from points $\{\ab_N\}\subset D_N$, that is $X^{(N)}_0 = \ab_N$ for $N\geq 1$.  We will assume that $\{\ab_N\}$ satisfies $\ab_N/N\rightarrow \ab$ for some $\ab\in D$.

Let now $\tau_N$ be the exit time from $D_N$ by the random walk $\{X^{(N)}_n\}$, that is
$$\tau_N = \inf\big\{ n\geq 0:  X^{(N)}_n \not\in D_N\big\}.$$  
In addition, for $r\geq 1$, let $T^{(r)}_{N,\X}$ be the $r$th hitting time of $\X\in \Z^d$.
When $r=1$, we will denote $T_{N,x} = T^{(1)}_{N,x}$.
Define also $\tau_{\ab, D}$ as the exit time from $D$ of a $d$-dimensional Brownian motion, starting from $\ab\in D$.

With this notation, the range of the random walk $X^{(N)}_\cdot$ up to the time of exit is given by
$$R_N \ = \ \mathcal{R}_{\tau_N} \ = \ \sum_{\X\in D_N} 1\big(T_{N,\X}<\tau_N\big)$$
and the associated $p$-multiple range is given by
$$R^{(p)}_N \ = \ \mathcal{R}^{(p)}_{\tau_N} \ = \ \sum_{\X\in D_N} 1\big(T^{(p)}_{N,\X}<\tau_N<T^{(p+1)}_{N,\X}\big).
$$

Finally, when $d\geq 3$, let $p_0<1$ be the probability that (transient) simple random walk on $\Z^d$, starting from the origin, returns.

\begin{theorem}[Range]
\label{thm1}
We have the weak convergence, as $N\uparrow\infty$, that
\begin{align*}
\frac{R_N}{N^2/\log N} &\Rightarrow \pi \tau_{\ab,D} \ \ {\rm when \ }d=2, \ \ {\rm and \ }\\
\frac{R_N}{N^2} &\Rightarrow \frac{d}{2}\big(1-p_0\big)\tau_{\ab, D} \ \ {\rm when \ }d\geq 3.
\end{align*}

Moreover, for $k\geq 1$, we have that the $k$th moments of $R_N/(N^2/\log N)$ when $d=2$ and those of $R_N/N^2$ when $d\geq 3$ converge to the $k$th moments of their distributional limits, which satisfy a system of Poisson PDE via \eqref{helms_eq} and \eqref{mom_lim}.

\end{theorem}

\begin{theorem}[Multiple range]
\label{thm2}
Let $p\geq 1$.  We have the weak convergence, as $N\uparrow\infty$, that
\begin{align*}
\frac{R^{(p)}_N}{N^2/\log^2 N} &\Rightarrow 2\pi^2 \tau_{\ab,D} \ \ {\rm when \ }d=2, \ \ {\rm and \ }\\
\frac{R^{(p)}_N}{N^2} &\Rightarrow \frac{d}{2}\big(1-p_0\big)^{2}(p_0)^{p-1}\tau_{\ab, D} \ \ {\rm when \ }d\geq 3.
\end{align*}

Moreover, for $k\geq 1$, we have that the $k$th moments of $R^{(p)}_N/(N^2/\log^2 N)$ when $d=2$ and those of $R^{(p)}_N/N^2$ when $d\geq 3$ converge to the $k$th moments of their distributional limits, which satisfy a system of Poisson PDE via \eqref{helms_eq} and \eqref{mom1_lim}.

\end{theorem}

Let $\P_{N,\b}$ and $\E_{N,\b}$ be the probability measure and expectation governing the random walk path $X^{(N)}_\cdot$ where $X^{(N)}_0=\b$.  Let also $\P_\b$ and $\E_\b$ be the process measure and expectation with respect to $d$-dimensional standard Brownian motion starting from $\b$.  

Define, 
for $k\geq 1$, that
$$u^{(k)}(\ab) = 
\E_\ab\Big[\big(\tau_{\ab, D}\big)^k\Big].
$$
Such moments are known to be `polyharmonic', that is they satisfy a following system of PDE \cite{Helms}:  Let $k\geq 1$.  Then, $u^{(k)}(\ab) = 0$ for $\ab\not\in D$ and, for $\ab\in D$, we have
\begin{align}
\Delta u^{(1)} &= -2 \nonumber \\
\Delta u^{(k+1)} &= -2(k+1)u^{(k)}.
\label{helms_eq}
\end{align}
In particular, the limits 
\begin{align}
\lim_{N\uparrow\infty} \E_{N,\ab_N}\Big[\Big(R_N/(N^2/\log N)\Big)^k\Big] &= \pi^k u^{(k)}(\ab)  \
 \ {\rm when \ }d=2,\nonumber \\
\lim_{N\uparrow\infty} \E_{N,\ab_N}\Big[\Big(R_N/N^2\Big)^k\Big] &= \Big(\frac{d}{2}(1-p_0)\Big)^k u^{(k)}(\ab) \ \ {\rm when \ }d\geq 3.
\label{mom_lim}
\end{align}

Also, for $p\geq 1$,
\begin{align}
\lim_{N\uparrow\infty} \E_{N,\ab_N}\Big[\Big(R^{(p)}_N/(N^2/\log N)\Big)^k\Big] &= (2\pi^{2})^k u^{(k)}(\ab)  \
 \ {\rm when \ }d=2,\nonumber \\
\lim_{N\uparrow\infty} \E_{N,\ab_N}\Big[\Big(R^{(p)}_N/N^2\Big)^k\Big] &= \Big(\frac{d}{2}(1-p_0)^{2}(p_0)^{p-1}\Big)^k u^{(k)}(\ab) \ \ {\rm when \ }d\geq 3.
\label{mom1_lim}
\end{align}

\begin{remark}[Generalizations]
\label{rmk1}
\rm
We comment on some generalizations.

{\it Mean-zero walks.}  The proofs for Theorems \ref{thm1} and \ref{thm2} carry over straightforwardly to finite-range mean-zero random walks, with covariance matrix $\Gamma$, but now using also $\{\cos(n\pi \cdot): n\geq 0\}$ in addition to $\{\sin(n\pi \cdot): n\geq 0\}$ as bases functions in the Fourier decompositions.  The statements would be similar to Theorems \ref{thm1} and \ref{thm2}, but the term $d\tau_{\ab, D}$ would be replaced by the exit time from $D$ of a Brownian motion with covariance matrix $\Gamma$ starting from $\ab$.

{\it Biased walks.}  One might consider biased finite-range random walks with (vector) mean $\mathbf{m}\neq 0$.  Adapting the proofs given here, one can show, in $d\geq 2$, that $R_N/N \rightarrow (1-p_0)c(\ab,D)$ and $R^{(p)}_N/N \rightarrow (1-p_0)^2 (p_0)^pc(\ab, D)$ in probability, where $c(\ab,D) = \inf \{c > 0\, |\, \ab + c \, \mathbf{m} \in D^c\}$. 
\end{remark}

\begin{remark}[Dimension $d=1$]
\label{1D}\rm
To compare, we specify the limit for the scaled range up to the time of exit from $(0,N)$ in $d=1$ among other local time results in \cite{locker}:  Namely $R_N/N\Rightarrow \zeta_{\ab}$ where $\zeta_{\ab}\neq \tau_{\ab, (0,1)}$ and has density $f(x)=(\ab\wedge (1-\ab))/x^2$ for $\ab\wedge (1-\ab)<x<\ab\vee(1-\ab)$, $f(x) = 1/x^2$ for $\ab\vee(1-\ab)\leq x\leq 1$ and $f(x)=0$ for $x$ otherwise.

\end{remark}

In the next Section \ref{proofs}, the proofs of Theorems \ref{thm1} and \ref{thm2} are given, with the aid of estimates in Sections \ref{line proofs} and \ref{p_moment_section}.

\section{Proof of Theorems \ref{thm1} and \ref{thm2}}
\label{proofs}

We now detail three ingredients, two of which are known, used in the proof of Theorems \ref{thm1} and \ref{thm2}.

\medskip
{\bf (A)  Asymptotics of $\mathcal{R}_n$ and $\mathcal{R}^{(p)}_n$}.  It is known when $d=2$ that $R_n/(n/\log n) \rightarrow \pi$ a.s. \cite{Dvoretzky}.  Whereas when $d\geq 3$, since the random walk is transient, $R_n/n\rightarrow 1-p_0$ a.s. where $p_0$ is the probability of return to the starting point \cite{Dvoretzky}, \cite[p. 38-40]{Kesten}.  

Also, for $p\geq 1$, when $d=2$, it is known that $\mathcal{R}^{(p)}_n/(n/\log^2(n)) \rightarrow \pi^2$ a.s. \cite{Flatto}.  But, when $d\geq 3$, $\mathcal{R}^{(p)}_n/n \rightarrow (1-p_0)^{2}(p_0)^{p-1}$ a.s. \cite{Pitt}. 

We note all of these results do not depend on the starting point value as the random walk dynamics on $\Z^d$ is translation-invariant.
\medskip

{\bf (B)  Input from a functional CLT}.  
When $\ab_N/N \rightarrow \ab\in D$, we have by the functional central limit theorem that the random walk paths $\{\frac{\sqrt{d}}{N}X^{(N)}_{[N^2s]}: s\geq 0; X^{(N)}_0=\ab_N\}$ converge weakly say in the uniform topology to Brownian motion $\{B_s: s\geq 0; B_0=\ab\}$ (cf. Sections 16, 18 \cite{Billingsley}).  

Since $D$ is an open, bounded and Lipschitz (or Jordan in $d=2$) domain, the time $\tau_{\ab, D}<\infty$ is a continuous function with respect to the uniform topology on the space of Brownian trajectories a.s.  Indeed, let $\omega$ be a Brownian trajectory starting from $\ab$, and $\{\omega^n\}$ be a sequence of continuous paths converging to it uniformly on compact time intervals.  One cannot have the limit $u=\lim \tau_{\ab, D}(\omega^n)< \tau_{\ab, D}(\omega) = v$, since then, from the uniform convergence, $u \geq \tau_{\ab, D}(\omega) = v$.  But, one cannot have $u>v$ either, since from the uniform convergence, $\omega(r)\in \bar{D}$ for $v\leq r\leq u$, a contradiction given that $\omega$ must also visit $\R^d\setminus \bar{D}$ in this time interval as the Lipschitz (or Jordan in $d=2$) boundary $\partial D$ satisfies (after a conformal transformation in $d=2$) a uniform cone condition (cf.  \cite[Ch. 4]{A}, \cite[Ch. 4]{KS}).

   Then, by the continuous mapping theorem, we have that
\begin{equation}
\label{fclt}
\frac{\tau_N}{N^2} \Rightarrow d\tau_{\ab,D}.
\end{equation}
 In particular, we have the convergence in probability as an immediate consequence,
 \begin{equation}
 \label{log}
 \frac{\log \tau_N}{\log N} \stackrel{P}{\longrightarrow} 2.
 \end{equation}

\medskip

{\bf (C)  Moment estimates}.  We show in Sections \ref{line proofs} and \ref{p_moment_section}, with the aid of a `conductance' estimate, the following limits.  When $\ab_N/N\rightarrow \ab\in D$,
we claim that
\begin{align}
\label{line1}
 \limsup_{N\uparrow \infty} \E_{N,\ab_N} \Big[\Big(\frac{R_N}{N^2/\log N}\Big)^k\Big] <\infty & \quad {\rm when \ }d=2\\
\label{line2}
\limsup_{N\uparrow\infty} \E_{N,\ab_N}\Big[\Big(\frac{R_N}{N^2}\Big)^k\Big] <\infty& \quad {\rm when \ }d\geq 3.
\end{align}

Moreover, we claim that
\begin{align}
\label{linep1}
\limsup_{N\uparrow\infty} \E_{N,\ab_N} \Big[\Big(\frac{R^{(p)}_N}{N^2/\log^2 N}\Big)^k\Big] <\infty & \quad {\rm when \ }d=2\\
\label{linep2}
\limsup_{N\uparrow\infty} \E_{N,\ab_N}\Big[\Big(\frac{R^{(p)}_N}{N^2}\Big)^k\Big] <\infty& \quad {\rm when \ }d\geq 3.
\end{align}

Although not needed for Theorems \ref{thm1} and \ref{thm2}, corresponding positive lower bounds can also be shown by arguments similar to those given for the upper bounds:
\begin{equation}\liminf_{N\uparrow\infty}\E_{N,\ab_N}[(R_N/w_N)^k] >0 \ \ {\rm and  \ \ } \liminf_{N\uparrow\infty}\E_{N,\ab_N}[(R^{(p)}_N/v_N)^k]>0,
\label{eq:lower_bounds}
\end{equation}
 where $w_N = N^2/\log(N)$ and $v_N = N^2/\log^2(N)$ in $d=2$, and $w_N=v_N =N^2$ in $d\geq 3$.    As a note to the interested reader, by Jensen's inequality, it suffices to show these claims \eqref{eq:lower_bounds} for $k=1$.

\medskip
{\bf Proof of Theorems \ref{thm1} and \ref{thm2}.} In terms of the `ingredients', we now combine them in the following way.
Suppose $d=2$ and write
$$
\frac{R_N}{N^2/\log N} = \frac{\mathcal{R}_{\tau_N}}{\tau_N/\log \tau_N} \frac{\tau_N}{N^2}\frac{\log N}{ \log \tau_N}.
$$
since $\tau_N \uparrow \infty$ a.s., by ingredient {\bf A}, we have that
${\mathcal{R}_{\tau_N} \log \tau_N}/{\tau_N} \rightarrow \pi$ converges in probability.  On the other hand, from ingredient {\bf B}, we have that $(\tau_N/N^2)(\log N/\log \tau_N)\Rightarrow \tau_{\ab,D}$.  Hence,
$R_N/(N^2/\log N) \Rightarrow \pi \tau_{\ab, D}$ as desired.

The weak convergence argument for $R_N/N^2 \Rightarrow (d/2)(1-p_0)\tau_{\ab, D}$ when $d\geq 3$ is similar.

Then, with respect to Theorem \ref{thm1}, convergence of the moments follows immediately from the weak convergence and ingredient {\bf C}. 

Finally, we note that the argument for Theorem \ref{thm2} follows the same steps.  \qed

     \section{Moment bounds: Proofs of \eqref{line1} and \eqref{line2}}
     \label{line proofs}
We will prove only \eqref{line1}, as the argument for \eqref{line2} is easier, and in particular can be done in a similar way.		Fix $d=2$ for the remainder of the section.

    The bound \eqref{line1} will hold if we establish
    \begin{align}
    \label{line1.1}
    B:= \limsup_{N\uparrow\infty} \sup_{\b\in D_N} \frac{\log N}{N^2} \E_{N,\b}\big [R_N \big] < \infty
    \end{align}
    and, for $k\geq 2$, the factorial moment
    \begin{align}
    \label{line1.2}
     \limsup_{N\uparrow\infty}\sup_{\b\in D_N} \Big(\frac{\log(N)}{N^2}\Big)^k\E_{N,\b}\Big[R_N (R_N-1)\cdots (R_N- (k-1))\Big] \leq B^k<\infty.
     \end{align}

\noindent {\bf Preliminaries.} Consider now the random walk hitting probability $P^{(N)}_{\b}: \Z^2 \rightarrow [0,1]$, for $\b\in D_N$, defined by  
\begin{equation}
\label{p_def}
P^{(N)}_{\b}(\X) = \P_{N,\b}(T_{N,\X}<\tau_N)
\end{equation}
 where $T_{N,\X}$ is the random walk hitting time of $\X\in \Z^2$.

 We now make a reduction step: Since $D\subset (-\frac{A}{2}, \frac{A}{2})^d$, for an integer $A<\infty$, we may bound 
     $\P^{(N)}_\b(\X) \leq \P_{N,\b}(T_{N,\X}<\tau'_N)$ where $\tau'_N$ is the exit time from $NA(-\frac{1}{2},\frac{1}{2})^2$.  Then,
         $R_N \leq \mathcal{R}_{\tau'_{NA}}$ where $\tau'_{NA}$ is the time to exit from a cube of size 1, scaled by $NA$.  Hence, it will be enough to show the bounds \eqref{line1}, \eqref{line2}, via translation-invariance, with respect to the cube $(0,1)^2$.  Accordingly, we fix for the remainder of the section that $D= (0,1)^2$.

We now derive an explicit expression for $P^{(N)}_\b$.  Note $P^{(N)}_{\b} (\b) = 1$ and, for $\X\not\in D_N$, that 
$P^{(N)}_{\b} (\X) = 0$.
On the other hand, if $\X\in D_N$ and $\b\neq \X$, then first step analysis yields
\begin{equation}
\label{p_harm}
P^{(N)}_{\b}(\X) = \frac{1}{4} \left( P^{(N)}_{\b+ \e_1}(\X) + P^{(N)}_{\b-\e_1}(\X) + P^{(N)}_{\b+\e_2}(\X) + P^{(N)}_{\b-\e_2}(\X) \right).
\end{equation}

Define a notion of `conductance' $g_N : \,  \Z^2\to [0,1]$ by
\begin{align}
g_N (\X)  &= 1 - \frac{1}{4} \left( P^{(N)}_{\X+\e_1}(\X) + P^{(N)}_{\X-\e_1}(\X) + P^{(N)}_{\X+\e_2}(\X) + P^{(N)}_{\X-\e_2}(\X) \right),\nonumber \\
&= \P_\X\big(\tau_N<T_{N,\X}\big),
\label{g_prob}
\end{align}
so that for $\X \in D_N$, we have
\begin{equation}
\label{one}
P^{(N)}_{\b}(\X) - \frac{1}{4} \left( P^{(N)}_{\b + \e_1}(\X) + P^{(N)}_{\b-\e_1}(\X) + P^{(N)}_{\b+\e_2}(\X) + P^{(N)}_{\b-\e_2}(\X) \right) = g_N (\X) \, 1_{\b}(\X). 
\end{equation}

Write $\b = (a,b)$ and consider now the discrete Fourier transforms of each term in \eqref{one}:
\begin{align*}
P^{(N)}_{\b} (\X) &= \sum_{m=1}^{N-1}\sum_{n=1}^{N-1} a_{mn} \sin\left(\frac{n \pi a}{N}\right)\sin\left(\frac{m \pi b}{N}\right),
\end{align*}
\begin{align*}
P^{(N)}_{\b \pm \e_1} (\X) &= \sum_{m=1}^{N-1}\sum_{n=1}^{N-1} a_{mn} \sin\left(\frac{n \pi (a\pm 1)}{N}\right)\sin\left(\frac{m \pi b}{N}\right),\\
P^{(N)}_{\b \pm \e_2} (\X) &= \sum_{m=1}^{N-1}\sum_{n=1}^{N-1} a_{mn} \sin\left(\frac{n \pi a}{N}\right)\sin\left(\frac{m \pi (b\pm 1)}{N}\right),\\
g_N (\X) 1_{\b}(\X) &= \sum_{m=1}^{N-1}\sum_{n=1}^{N-1} b_{mn} \sin\left(\frac{n \pi a}{N}\right)\sin\left(\frac{m \pi b}{N}\right).
\end{align*}
Through straightforward trigonometric manipulations, \eqref{one} is re-expressed as
\begin{align}
\label{two}
\sum_{m=1}^{N-1}\sum_{n=1}^{N-1} a_{mn} \sin\left(\frac{n \pi a}{N}\right)\sin\left(\frac{m \pi b}{N}\right) \left(1-\frac{1}{2}\left[\cos\left(\frac{n \pi}{N}\right) + \cos\left(\frac{m \pi}{N}\right) \right] \right)= g_N (\X) \, 1_{\b}(\X).
\end{align}

We now find formulas for $a_{mn}$ and $b_{mn}$.  Recall the orthogonality relation, for $n,n' \in \N$, that
$$\sum_{\ell=1}^{N-1} \sin\left(\frac{n \pi \ell}{N}\right)\sin\left(\frac{n' \pi \ell}{N}\right) = \frac{N1(n=n')}{2}.$$
Then, 
$$\sum_{r=1}^{N-1}\sum_{\ell=1}^{N-1} g_N (\X) \, 1_{\b}(\X) \, \sin\left(\frac{n' \pi \ell}{N}\right)\sin\left(\frac{m' \pi r}{N}\right) = \frac{N^2b_{m'n'}}{4}.$$
Hence, writing $\X = (x_1, x_2)$, we have
$$b_{mn} = \frac{4 g_N (x,y)\sin\left(\frac{n \pi x_1}{N}\right)\sin\left(\frac{m \pi x_2}{N}\right)}{N^2}.$$
Moreover, by equating the coefficients $a_{mn}[(1-(1/2)[\cos(n\pi/N) + \cos(m\pi/N)])] = b_{mn}$ with respect to \eqref{two}, we have
$$a_{mn} = \frac{4g_N (\X)\sin\left(\frac{n \pi x_1}{N}\right)\sin\left(\frac{m \pi x_2}{N}\right)}{N^2\left(1-\frac{1}{2}\left[\cos\left(\frac{n \pi}{N}\right) + \cos\left(\frac{m \pi}{N}\right) \right] \right)}.$$

We now obtain explicit formulas for $g_N(\X)$ and $P^{(N)}_{\b}$.  The boundary condition $P^{(N)}_\X(\X) = 1$ yields
$$P^{(N)}_{\X} (\X) = \sum_{m=1}^{N-1}\sum_{n=1}^{N-1} a_{mn} \sin\left(\frac{n \pi x_1}{N}\right)\sin\left(\frac{m \pi x_2}{N}\right) = 1,$$
from which
$$g_N (\X) = \left[\sum_{m=1}^{N-1}\sum_{n=1}^{N-1} \frac{4\sin^2\left(\frac{n \pi x_1}{N}\right)\sin^2\left(\frac{m \pi x_2}{N}\right)}{N^2\left(1-\frac{1}{2}\left[\cos\left(\frac{n \pi}{N}\right) + \cos\left(\frac{m \pi}{N}\right) \right] \right)}\right]^{-1}
 $$
and
\begin{equation}
\label{P_N formula}
P^{(N)}_{\b}(\X) = \sum_{m=1}^{N-1}\sum_{n=1}^{N-1}\frac{4g_N (\X)\sin\left(\frac{n \pi x_1}{N}\right)\sin\left(\frac{m \pi x_2}{N}\right)\sin\left(\frac{n \pi a}{N}\right)\sin\left(\frac{m \pi b}{N}\right)}{N^2\left(1-\frac{1}{2}\left[\cos\left(\frac{n \pi}{N}\right) + \cos\left(\frac{m \pi}{N}\right) \right] \right)}.
\end{equation}

We now claim the following upper bound, proved at the end of the subsection,
\begin{equation}\label{g_N}
\limsup_{N\uparrow\infty}\sup_{N/\log^2(N) \leq x_1,x_2\leq N- N/\log^2(N)} \log(N) g_N(\X) < \infty.
\end{equation}
\vskip .1cm

We are now in position to show \eqref{line1.1} and \eqref{line1.2}.

\noindent {\bf Proof of \eqref{line1.1}.}
We may write the expected range $\E_{N,\b}[R_N]$ as
$$\E_{N,\b}[R_N] = \E_{N,\b}\left[ \sum_{\X \in D_N} 1(T_{N,\X} < \tau_N)\right] = \sum_{\X\in D_N} P^{(N)}_{\b}(\X).$$
Let 
\begin{equation}
\label{D_1_set}
D_N^1 = \{\X\in D_N:  N/\log^2(N) \leq x_1,x_2 \leq N - N/\log^2(N)\}.
\end{equation}
 Then, 
\begin{align*}
\E_{N,\b}[R_N] 
&= \sum_{ \X\in D^1_N} P^{(N)}_\b(\X) + \sum_{\X\in D_N\cap (D_N^1)^c}P^{(N)}_\b(\X) = J_1+J_2.
\end{align*}
The second sum $J_2$, bounding $P^{(N)}_\b(\X)\leq 1$, is of order $O(N^2/\log^2(N))$.  

We now show an $O(N^2/\log(N))$ bound for the first sum $J_1$ to finish.
Note the limit
$$\lim_{N\rightarrow\infty}
\frac{1}{N}\sum_{\ell=1}^N \sin(r\pi \ell /N) = \int_0^1 \sin(r\pi u)du= \frac{1(r \ {\rm odd})}{\pi r}.$$
With the formula \eqref{P_N formula} and bound \eqref{g_N} in hand, uniformly over $\b$, we conclude
\begin{align*}
J_1&\leq O(N^2/\log N) \sum_{\X\in D^1_N} \sum_{1\leq n,m\leq N} \frac{\sin(n\pi x_1/N)\sin(m\pi x_1/N)}{N^2[n^2 + m^2]}\\
& \leq O(N^2/\log N) \sum_{1\leq n,m\leq N} \frac{1}{nm(n^2+ m^2)} = O(N^2/\log N),
\end{align*}
as desired. \qed
\vskip .2cm

\noindent {\bf Proof of \eqref{line1.2}.}
For $1\leq \ell \leq k-1$, define and bound
$$P^{(N)}_{\X}(\y; \{\z_\ell, \ldots, \z_k\}) = \P_\X ( T_{N,\y} < \min\{ T_{N,\z_\ell}, \ldots, T_{N,\z_k}, \tau_N\}) \leq P^{(N)}_\X(\y).$$
Denote by $\S_k$ the set of permutations of $\{1,2,\ldots, k\}$.  The factorial moment of $R_N$ of order $k$ is writtten as
\begin{align*}
&\E_{N,\b}\Big[R_N\big(R_N-1\big)\cdots \big(R_N-(k-1)\big)\Big] \\
&\ \ = \E_{N,\b}\Big[\mathop{\sum_{\z_1,\ldots, \z_k\in D_N}}_{{\rm distinct}}  \prod_{\ell=1}^k 1(T_{\z_\ell}<\tau_N)\Big]\\
&\ \ =\mathop{\sum_{\z_1,\ldots, \z_k\in D_N}}_{{\rm distinct}} \sum_{\pi \in \S_k}  \P_\b(\z_{\pi_1}; \{\z_{\pi_2}, \ldots, \z_{\pi_k}\})\P_{\z_{\pi_1}}(\z_{\pi_2}; \{\z_{\pi_3},\ldots, \z_{\pi_k}\})\cdots \P_{\z_{\pi_{k-1}}}(\z_{\pi_k}).
\end{align*}

We now observe, for each $\pi\in \S_k$, by the proven \eqref{line1.1}, that
\begin{align*}
\mathop{\sum_{\z_1,\ldots, \z_k\in D_N}}_{{\rm distinct}}
& \P_\b(\z_{\pi_1}; \{\z_{\pi_2}, \ldots, \z_{\pi_k}\})\P_{\z_{\pi_1}}(\z_{\pi_2}; \{\z_{\pi_3},\ldots, \z_{\pi_k}\})\cdots \P_{\z_{\pi_{k-1}}}(\z_{\pi_k})\\
&\ \ \leq 
\mathop{\sum_{\z_1,\ldots, \z_k\in D_N}}_{{\rm distinct}}
 \P_\b(\z_{\pi_1})\P_{\z_{\pi_1}}(\z_{\pi_2})\cdots \P_{\z_{\pi_{k-1}}}(\z_{\pi_k})\\
&\ \ = \E_{N,\b}[R_N]\prod_{\ell=2}^k \E_{N,\z_{\pi_{\ell-1}}}[R_N] = O\big( (N^2/\log(N))^k\big).
\end{align*}
Hence, as a consequence, \eqref{line1.2} follows immediately.
\qed
\vskip .2cm

\noindent {\bf Proof of \eqref{g_N}.}
Since $1-\cos(n\pi/N) = 2\sin^2(n\pi/(2N))$,
$2\sin^2(\ell\pi/(2N)) \leq 2\ell^2\pi^2/(4 N^2)$,
we have
\begin{align*}
g_N^{-1}(\X) &\geq 
\frac{2}{\pi^2}\sum_{n=1}^{N}\sum_{m=1}^{N} \frac{\sin^2(n\pi x_1/N)\sin^2(m\pi x_2/N)}{n^2 + m^2}.
\end{align*}
Note, when $\X$ is the midpoint, $\X=(\lfloor N/2\rfloor, \lfloor N/2\rfloor)$,
$$g_N^{-1}(\X) \ \geq \ \frac{2}{\pi^2}\mathop{\sum_{1\leq n,m\leq N}}_{{\rm odd}}\frac{1}{n^2 + m^2} \ \geq (1/\pi^2)\log(N).$$

When $\X$ is not the midpoint, let $d_N=d_N(\X)$ be its distance to the boundary of $D_N$.  Let $B_N=B_N(\X)\subset D_N$ be the cube with width $d_N$ and center $\X$.  Let $\zeta_N=\zeta_N(\X)$ be the time for the random walk to exit $B_N$.  Since $\zeta_N\leq \tau_N$, by \eqref{g_prob}, we have that
$g_N(\X) \leq \P_\X(\zeta_N<T_{N,\X}):= \widetilde g_N(\X)$.

One can derive, a formula and lower bound for $\widetilde g_N(\X)$, as for $g_N(\X)$ when $\X$ is the midpoint above, using the scale $d_N$ and domain $B_N$ instead of $N$ and $D_N$:
\begin{align*}
\widetilde g^{-1}_N(\X) &= \sum_{n=1}^{d_N}\sum_{m=1}^{d_N} \frac{4\sin^2\big(\frac{n\pi x_1}{d_N}\big)\sin^2\big(\frac{n\pi x_2}{d_N}\big)}{N^2\big(1-\frac{1}{2}\big[\cos\big(\frac{n\pi}{d_N}\big) + \cos\big(\frac{m\pi}{d_N}\big)\big]\big)}\\
&\geq    (1/\pi^2)\log(d_N).
\end{align*}

Hence, when $N/\log^2(N) \leq x_1,x_2  \leq N - N/\log^2(N)$, we have that
\begin{equation}
\label{tilde g_N}
g_N^{-1}(\X) \geq \widetilde g^{-1}_N(\X)\geq (1/2\pi^2)\log(N)
\end{equation}
 for all large $N$, as desired. \qed

\section{Moment bounds:  Proofs of \eqref{linep1} and \eqref{linep2}}
\label{p_moment_section}

For $d \geq 3$, the bound in~\eqref{linep2} follows from~\eqref{line2} since $R_N^{(p)} \leq R_N$ a.s.
So we only need to show~\eqref{linep1} when $d=2$. We write 
$$\E_{N,\b}[R^{(p)}_N] \ = \ \sum_{\X\in D_N} \E_{N,\b}\Big[ 1\big(T^{(p)}_{N,\X}<\tau_N<T^{(p+1)}_{N,\X}\big)\Big]
\ = \ \sum_{\X\in D_N} \Q_\b^{(N)}(\X)$$
where $\Q_\b^{(N)}:\Z^d\rightarrow [0,1]$ is given by 
\begin{align*}
\Q^{(N)}_\b(\X) &= \P_{N,\b}(T^{(p)}_{N,\X}<\tau_N<T^{(p+1)}_{N,\X})\\
&\leq \P_{N,\b}(T^{(p)}_{N,\X}<\tau_N)\P_{N,\X}(\tau_N<T_{N,\X}).
\end{align*}
As $P^{(N)}_\b(T^{(p)}_{N,\X}<\tau_N) \leq \P_{N,\b}(T_{N,\X}<\tau_N)=P^{(N)}_\b(\X)$ and $g_N(\X)=\P_{N,\X}(\tau_N<T_{N,\X})$ (cf. \eqref{g_prob}), we have
$$\Q^{(N)}_\b(\X) \leq g_N(\X)\P^{(N)}_\b(\X).$$

Recall $\widetilde g_N$ from the proof of \eqref{g_N}.  Note that $g_N(\X)\leq \widetilde g_N(\X)$ and that $\widetilde g_N(\X) =O(\log^{-1}(N))$ when $d_N$, the distance between $\X$ and the boundary of $D_N$, is greater than $N/\log^2(N)$ say. 

Consider now the set $D^2_N\subset D_N$ of points $\X$ away by at least $N/\log^2(N)$ from the boundary of $D_N$. 
From our assumptions on $D$, since $\partial D$ has finite perimeter, the area of the region within $N\log^{-2}(N)$ of $\partial D$ is of order $O({\rm Per}(\partial D)\, N \log^{-2}(N))$, and so the number of points in $D_N$ within distance $N\log^{-2}(N)$ is of order $O(N^2/\log^{2}(N))$.

Recall, by the proven \eqref{line1} that $\E_{N,\b}[R_N]=O(N^2/\log(N))$ uniformly over $\b\in D_N$.
We then obtain, uniformly over $\b\in D_N$, that
\begin{align*}
\E_{N,\b}[R^{(p)}_N] &\leq O(N^2/\log^2(N)) + \sum_{\X\in D^2_N} \widetilde g_N(\X)P^{(N)}_\b(\X)\nonumber \\
&\leq O(N^2/\log^2(N)) + O(\log^{-1}(N))\E_{N,\b}[R_N] \  = O(N^2/\log^2(N)).\nonumber
\end{align*}
This gives \eqref{linep1} with respect to $k=1$.

We now turn to estimating the factorial moment of order $k\geq 2$.  We concentrate on the case $k=2$ which encapsulates the main ideas, and then comment on the case $k\geq 3$.

Consider that 
\begin{align}
&\P_{N,b}(T^{(p)}_{N,\X}, T^{(p)}_{N,\y}<\tau_N< T^{(p+1)}_{N,\X}, T^{(p+1)}_{N,\y}) \nonumber\\
&\ \ \ \ \ = \P_{N,b}(T^{(p)}_{N,\X}< T^{(p)}_{N,\y}<\tau_N< T^{(p+1)}_{N,\X}, T^{(p+1)}_{N,\y}) +\P_{N,b}(T^{(p)}_{N,\y}< T^{(p)}_{N,\X}<\tau_N< T^{(p+1)}_{N,\X}, T^{(p+1)}_{N,\y}).
\label{multirangehelp0}
\end{align}

We bound now the first term on the right-hand side of \eqref{multirangehelp0}.  Let again $B_N$ be a cube centered at $x$ with width $d_N(\X)=O(N/\log^2(N))$ and $\zeta_N$ be the exit time from $B_N$.  In decomposing the path specified in the first term, we argue, for $y\not\in B_N$, that
\begin{align}
&\P_{N,b}(T^{(p)}_{N,\X}< T^{(p)}_{N,\y}<\tau_N< T^{(p+1)}_{N,\X}, T^{(p+1)}_{N,\y})\nonumber \\
&\ \ \ \ \ \leq p\times \P_{N,\b}(T_{N,\X}<\tau_N)\E_{N,\X}[1(\zeta_N<T_{N,\X}) \P_{N,X^{(N)}_{\zeta_N}}(T_{N,\y}<\tau_N)]
\P_{N,\y}(\tau_N<T_{N,\y}).
\label{multirangehelp1}
\end{align}
We obtain a similar expression bounding the second term in \eqref{multirangehelp0} by interchanging $\X$ and $\y$.
Indeed, the left-hand side of~\eqref{multirangehelp1} is the probability of the set of paths that start at $b$, hit $x$ exactly $p$ times before the $p$-th visit to $y$ after which the boundary $\partial D_N$ is hit before coming back to $x$ or $y$. The set of these paths is a subset of paths, starting from $b$, hitting $x$ (exactly $p$ times and hitting $y$ at most $\ell(\y)\leq p-1$ times), which then exit a box $B_N$ centered at $x$ (without hitting $x$ again) to hit $y$ (exactly $p-\ell(\y)$ times) before hitting $\partial D_N$. The probability of the set of such paths is bounded by the right-hand side expression.

By \eqref{tilde g_N}, $\P_{N,\X}(\zeta_N<T_{N,\X}) = O(\log^{-1}(d_N(\X)))$.  Summing \eqref{multirangehelp1}, and its analog with $\X$ and $\y$ interchanged, over $\X$ and $\y$, separated by at least $N/\log^2(N)$ and also when both $\X,\y$ are away from the boundary of $D_N$ by $N\log^2(N)$, noting the proved estimate \eqref{line1.1}, we obtain uniformly over $\b\in D_N$
the desired bound on the factorial moment
$$\E_{N,\b}[R^{(p)}_N(R^{(p)}_N-1)] = O(N^2/\log^2(N)) + O(\sup_{\z\in D_N}\E_{N,\z}[R_N]^2/ \log^{2}(N)) =
 O(N^4/\log^4(N)).$$

When $k\geq 3$, analogously, we may bound $k!$ terms such as 
\begin{equation}
\P_{N,\b}(T^{(p)}_{N,\z_1}<T^{(p)}_{N,\z_2}<\cdots <T^{(p)}_{N,\z_k}<\tau_N< T^{(p+1)}_{N,\z_1},\ldots, T^{(p+1)}_{N,\z_k}).
\label{multirangehelp2}
\end{equation}
Indeed, consider cubes of width $O(N/\log^2(N))$ centered around $\z_1, \ldots, \z_{k-1}$.  Let $\zeta_{1,N}, \ldots, \zeta_{k-1,N}$ be the exit times from these cubes.  Then, for $\z_1,\ldots, \z_k$ separated from each other and the boundary of $D_N$ by $O(N/\log^2(N))$, we have \eqref{multirangehelp2} is bounded by
\begin{align*}
&p^{k-1}\times \P_{N,\b}(T_{N,\z_1}<\tau_N)\E_{N,\z_1}[1(\zeta_{1,N}<T_{N,\z_1}) \P_{N,X^{(N)}_{\zeta_{1,N}}}(T_{N,\z_2}<\tau_N)]\\
&\ \ \ \ \ \ \times \E_{N,\z_2}[1(\zeta_{2,N}<T_{N,\z_2})\P_{N,X^{(N)}_{\zeta_{2,N}}}(T_{N,\z_3}<\tau_N)]\cdots \\
&\ \ \ \ \ \ \cdots \times \E_{N,\z_{k-2}}[1(\zeta_{k-1,N}<T_{N,\z_{k-1}})\P_{N,X^{(N)}_{\zeta_{k-1,N}}}(T_{N,\z_k}<\tau_N)]
 \P_{N,\z_{k}}(\tau_N<T_{N,\z_k}).
\end{align*}
Through the estimate \eqref{tilde g_N} again, summing such an expression over $\z_1,\ldots, \z_k$ we can bound the factorial moment uniformly over $\b\in D_N$ as 
\begin{align*}
\E_{N,\b}[R^{(p)}_N(R^{(p)}_N-1)\cdots (R^{(p)}_N - (k-1))] &= O(N^k/\log^2(N)) + O((\sup_{\z\in D_N}\E_{N,\z}[R_N])^k/\log^k(N)) \\
&= O((N^2/\log^2(k))^k).
\end{align*}

Hence, in this way, \eqref{linep1} follows.
\qed

\medskip

\noindent {\bf Acknowledgements.} 
We would like to thank Davar Khoshnevisan for useful conversations.  Thanks also to Antoine Jego for bringing to our attention references \cite{J1}, \cite{J2} after an initial version of this note was written.
SS was partially supported by ARO W911NF-18- 1-0311 and a Simons Foundation Sabbatical grant. SCV was partially supported by the Simons Foundation through awards 524875 and 560103.

 \bibliographystyle{plain}

\begin{thebibliography}{99}
 \frenchspacing
 
 \bibitem{AB}
 {\sc Abe, Y., Biskup, M.} (2019).
 \newblock Exceptional points of two-dimensional random walks at multiples of the cover time.
 \newblock arXiv:1903.04045
 
 \bibitem{A}
 {\sc Adams, R.A., Fournier, J.J.F.} (2003).
 \newblock {\it Sobolev Spaces.} 2nd Ed.
 \newblock Academic Press, Amsterdam.
 
 \bibitem{Asselah}
 {\sc Asselah, A., Schapira, B., Sousi, P.} (2019).
 \newblock Capacity of the range of random walk on $\Z^4$.
 \newblock {\it Ann. Probab.} {\bf 47} 1447--1497.
  
  \bibitem{locker}
 {\sc Athreya, S., Sethuraman, S., Toth, B.} (2011).
 \newblock On the range, local times and periodicity of random walk on an interval.
 \newblock {\it ALEA, Lat. Am. J. Probab. Math. Stat.} {\bf 8} 269--284.
 
 \bibitem{Bass}
 {\sc Bass, R.F., Chen, X., Rosen, J.} (2009).
 \newblock Moderate deviations for the range of planar random walks. 
 \newblock {\it Mem. Amer. Math. Soc.} {\bf 198} (929), viii+82.

\bibitem{Billingsley}
{\sc Billingsley, P.} (1968).
\newblock {\it Convergence of Probability Measures.}
\newblock Wiley, New York.
 
 \bibitem{Dvoretzky}
 {\sc Dvoretzky, A., Erd\"os, P.} (1950).
 \newblock Some problems on random walk in space; in {\it Proceedings of the Second Berkeley Symposium on Mathematical Statistics and Probability}, 353--367.
 \newblock University of California Press, Berkeley and Los Angeles.
 
\bibitem{Flatto}
{\sc Flatto, L.} (1976).
\newblock The multiple range of two-dimensional recurrent walk.
\newblock {\it Ann. Probab.} {\bf 4} 229--248.

\bibitem{Hamana}
{\sc Hamana, Y.} (1998).
\newblock A remark on the multiple point range of two dimensional random walk.
\newblock {\it Kyushu J. Math.} {\bf 52} 23--80.
 
 \bibitem{Helms}
 {\sc Helms, L.} (1967).
 \newblock Biharmonic functions and Brownian motion.
 \newblock {\it J. Appl. Probab.} {\bf 4} 130--136.
 
  

 
 \bibitem{Weiss}
 {\sc den Hollander, G.H., Weiss, F.} (1994).
 \newblock Trapping in transport processes; in {\it Contemporary problems in statistical physics}, 147--203. 
 \newblock Society for Industrial
and Applied Mathematics, Philadelphia.

\bibitem{J1} {\sc Jego, A.} (2019).
\newblock Characterization of planar Brownian multiplicative chaos.
\newblock arXiv:1909.05067v2

\bibitem{J2}
{\sc Jego, A.} (2020).
\newblock Thick points of random walk and the Gaussian free field.
\newblock {\it Elec. J. Probab.} {\bf 25} 1--39.

\bibitem{KS}
{\sc Karatzas, I., Shreve, S.} (1991).
\newblock {\it Brownian Motion and Stochastic Calculus}, 2nd Ed.
\newblock Springer-Verlag, New York.

 \bibitem{Pitt}
{\sc Pitt. J.H.} (1974).
\newblock Multiple points of transient random walks.
\newblock {\it Proc. Amer. Math. Soc.} {\bf 43} 195--199.

\bibitem{Kesten}
 {\sc Spitzer, F.} (1964).
 \newblock {\it Principles of random walk.}
 \newblock University Series in Higher Math, Van Nostrand, Princeton, NJ.





\end{thebibliography}
 
\end{document}